\documentclass[a4paper, leqno, 12pt]{amsart}
\usepackage{amssymb,amsthm,amsmath} 
\usepackage[T1]{fontenc} 
\usepackage[utf8]{inputenc}
\usepackage{color}
\usepackage[dvipsnames]{xcolor}
\usepackage[colorlinks, allcolors=RedViolet,pdfstartview=,pdfpagemode=UseNone]{hyperref}

\theoremstyle{plain}
\newtheorem{thm}[equation]{Theorem}
\newtheorem{lem}[equation]{Lemma}

\newtheorem{cor}[equation]{Corollary}

\numberwithin{equation}{section}
\theoremstyle{definition}
\newtheorem{defn}[equation]{Definition}
\newtheorem{eg}[equation]{Example}

\theoremstyle{remark}
\newtheorem{rem}[equation]{Remark}

\newcommand{\R}{\mathbb{R}}
\newcommand{\Phiw}{\Phi_{\text{\rm w}}}
\newcommand{\Phic}{\Phi_{\text{\rm c}}}

\newcommand{\essinf}{\operatornamewithlimits{ess\,inf}}
\newcommand{\supp}{\operatornamewithlimits{supp}}
\newcommand{\dist}{\operatorname{dist}}
\newcommand{\rad}{\operatorname{rad}}

\renewcommand{\epsilon}{\varepsilon}
\renewcommand{\phi}{\varphi}
\renewcommand{\rho}{\varrho}

\newcommand{\inc}[1]{\hyperref[def:aInc]{{\normalfont(Inc){\ensuremath{_{#1}}}}}}
\newcommand{\dec}[1]{\hyperref[def:aDec]{{\normalfont(Dec){\ensuremath{_{#1}}}}}}
\newcommand{\ainc}[1]{\hyperref[def:aInc]{{\normalfont(aInc){\ensuremath{_{#1}}}}}}
\newcommand{\adec}[1]{\hyperref[def:aDec]{{\normalfont(aDec){\ensuremath{_{#1}}}}}}
\newcommand{\azero}{\hyperref[def:A0]{{\normalfont(A0)}}}
\newcommand{\aone}{\hyperref[def:A1]{{\normalfont(A1)}}}
\newcommand{\aomega}{\hyperref[def:AOmega]{{\normalfont(A1)\ensuremath{_\Omega}}}}
\newcommand{\atwo}{\hyperref[def:A2]{{\normalfont(A2)}}}

\title{Extension in generalized Orlicz--Sobolev spaces}
\author{Jonne Juusti}

\begin{document}


\begin{abstract}
We study the existence of an extension operator $\Lambda \colon W^{1,\phi}(\Omega)\to W^{1,\psi}(\R^n)$.
We assume that $\phi \in \Phiw(\Omega)$ has generalized Orlicz growth, $\psi \in \Phiw(\R^n)$ is an extension of $\phi$, and that $\Omega\subset\R^n$ is an $(\epsilon,\delta)$-domain.
Special cases include the classical constant exponent case, the Orlicz case, the variable exponent case, and the double phase case.
\end{abstract}

\maketitle

%%%%%%%%%%%%%%%%%%%%%%%%%%%%%%%%%%%%%%%%%%%%%%%%%%%%%%%%
%%%%%%%%%%%%%%%%%%%%%%%%%%%%%%%%%%%%%%%%%%%%%%%%%%%%%%%%

\section{Introduction}

We study the existence of an extension operator $\Lambda : W^{1,\phi}(\Omega)\to W^{1,\psi}(\R^n)$.
We assume that $\phi \in \Phiw(\Omega)$ has generalized Orlicz growth, $\psi \in \Phiw(\R^n)$ is an extension of $\phi$, and $\Omega\subset\R^n$ is an $(\epsilon,\delta)$-domain.
Special cases of Orlicz growth include the  classical constant exponent case $\phi(x,t) = t^p$, the Orlicz case $\phi(x,t) = \phi(t)$, the variable exponent case $\phi(x,t) = t^{p(x)}$, and the double phase case $\phi(x,t) = t^p+a(x)t^q$.
Generalized Orlicz and Orlicz--Sobolev spaces on $\R^n$ have recently been studied in \cite{BenK21, HarK22, KamZ22, OrtR21} for example.

The existence of an extension operator is known in many special cases.
In \cite{Cal61},  Cald\'eron proved the existence of an extension in the classical constant exponent case under the conditions that $1 < p < \infty$ and $\Omega$ is a Lipschitz domain.
Stein extended these results to include the cases $p = 1,\infty$ in \cite{Ste70}, and Jones generalized the results further to $(\epsilon,\delta)$-domains in \cite{Jon81}.
In \cite{Chu92} Chua extended Jones' result to weighted constant exponent Sobolev spaces.
For a proof of the existence of an extension operator in the variable exponent case, see \cite[Chapter~8.5]{DieHHR11}.
Extension in the Orlicz case in general measure spaces has been studied in \cite{HeiT10}.

Our results include the cases mentioned above, and also other cases, like the double phase case.
It is worth mentioning, that the double phase extension operator has actually already been used in the proof of \cite[Proposition~6.5]{ColS16}.
The proof however does not cite any sources that show that such an operator actually exists, and, to best of our knowledge, our results are new in the double phase case.

The approach in this paper follows the approach used in \cite{DieHHR11}.
This means that we use extrapolation and the result of Chua in \cite{Chu92}.
One problem that is not present in the constant exponent and Orlicz cases, is that the function $\phi \in \Phiw(\Omega)$ is not defined outside $\Omega$.
Since $\phi$ depends on $x$, it is not trivial how to extend it to a function defined everywhere on $\R^n$.
To solve this problem we build upon the work in \cite{HarH19}, which gives an extension of $\phi$ in the case that $\Omega$ is quasi-convex.
Since $(\epsilon,\delta)$-domains are not necessarily quasi-convex, we prove a new result that suffices for our purposes.
This result is based on the fact that the conditions \aone{} and \aomega{} are equivalent in $(\epsilon,\delta)$-domains.
We will also show that the conditions \aone{} and \aomega{} are not equivalent in general, which answers a question left open in \cite{HarH19}.

The main reason for studying extension operators is that we can get easy proofs for some results that would be more difficult to prove otherwise.
As an example example, using extensions and the density of $C^\infty(\R^n)$ in $W^{1,\psi}(\R^n)$, we can prove that $C^\infty(\overline{\Omega})$ is dense in $W^{1,\phi}(\Omega)$.
The structure of the paper is as follows:
Section \ref{sec:preliminary} covers the basic properties of  generalized Orlicz and Orlicz--Sobolev spaces.
In section \ref{sec:OrliczExtension}, we study extending the domain of $\phi$, and also study the conditions \aone{} and \aomega{}.
Finally, the main result, that there exists an extension operator $\Lambda : W^{1,\phi}(\Omega)\to W^{1,\psi}(\R^n)$, is proved in Section \ref{sec:SobolevExtension}.

\section{Preliminaries}\label{sec:preliminary}

We briefly introduce our assumptions. More information about $L^\phi$-spaces can be found in \cite{HarH19}. 
A function $f : A \to [-\infty,\infty]$, where $A\subset\R$, is \emph{$L$-almost increasing}, if $f(s) \leq Lf(t)$ whenever $s < t$.
Similarly, $f$ is \emph{$L$-almost decreasing}, if $f(t)\leq Lf(s)$ whenever $s < t$.
By $\Omega$ we always denote an open subset of $\R^n$. 
If $f,g \colon \Omega\to[-\infty,\infty]$ are such that $f(x) \leq Cg(x)$  for some $C \geq 1$ and a.e $x \in \Omega$, then we write $f\lesssim g$. If $f \lesssim g \lesssim f$, then we write $f \approx g$.
We denote by $L^0(\Omega)$ the set of measurable functions in $\Omega$.

\begin{defn}
We say that $\phi : \Omega\times [0, \infty) \to [0, \infty]$ is a 
\emph{weak $\Phi$-function}, and write $\phi \in \Phiw(\Omega)$, if 
the following conditions hold:
\begin{itemize}
\item 
For every $x \in \Omega$ the function $t \mapsto \phi(x, t)$ is non-decreasing and
\[
\lim_{t \to 0^+}\phi(x,t) 
= 0
= \phi(x, 0)
\quad\text{and}\quad
\lim_{t \to \infty}\phi(x,t)
= \infty.
\]
\item 
There exists $L \geq 1$ such that $t \mapsto \frac{\phi(x,t)}{t}$ is $L$-almost increasing on $(0,\infty)$ for every $x \in \Omega$.
\item 
For every $f \in L^0(\Omega)$ the function $x \mapsto \phi(x,|f(x)|)$ is measurable.
\end{itemize}
If $t \mapsto \phi(x,t)$ is additionally convex and left-continuous for every $x \in \Omega$, then $\phi$ is a 
\emph{convex $\Phi$-function}, and we write $\phi \in \Phic(\Omega)$. If $\phi$ does not depend on $x$, then we omit the set  and write $\phi \in \Phiw$ or $\phi \in \Phic$.
\end{defn}

Two functions $\phi$ and $\psi$ are \emph{equivalent}, 
$\phi \simeq \psi$, if there exists $L \ge 1$ such that 
\[
\psi(x,L^{-1}t)
\leq \phi(x,t)
\leq \psi(x,Lt)
\]
for every $x \in \Omega$ and every $t \geq 0$.
Equivalent $\Phi$-functions give rise to the same space with 
comparable norms. 

We define the left-inverse of $\phi$ by
\[
\phi^{-1}(x,\tau)
:= \inf\lbrace t \geq 0 \colon \phi(x,t) \geq \tau \rbrace.
\]
Let us state some conditions for later use.

\begin{itemize}
\item[(A0)]\label{def:A0}
There exists $\beta \in (0, 1)$ such that
\[
\beta
\leq \phi^{-1}(x, 1)
\leq \beta^{-1}
\]
for every $x \in \Omega$.
\item[(A1)]\label{def:A1}
There exists $\beta \in (0,1)$ such that
\[
\beta \phi^{-1}(x,t)
\leq \phi^{-1}(y, t) 
\quad\text{when}\quad 
t
\in \left[1,\frac{1}{|B|}\right]
\]
for every ball $B$ with $|B| \leq 1$ and every $x,y\in B \cap \Omega$.
\item[(A2)]\label{def:A2}
For every $s > 0$ there exist $\beta \in (0, 1]$ and $h \in L^1(\Omega) \cap L^\infty(\Omega)$ such that
\[
\beta\phi^{-1}(x,t)
\leq \phi^{-1}(y,t)
\]
for every $x,y \in \Omega$ and every $t \in [h(x)+h(y),s]$.
\item[(aInc)$_p$] \label{def:aInc}
There exists $L_p \geq 1$ such that
\[
t
\mapsto \frac{\phi(x,t)}{t^{p}}
\]
is $L_p$-almost increasing in $(0,\infty)$ for every $x \in \Omega$.
\item[(aDec)$_q$] \label{def:aDec}
There exists $L_q \geq 1$ such that
\[
t
\mapsto \frac{\phi(x,t)}{t^{q}}
\]
is $L_q$-almost decreasing in $(0,\infty)$ for every $x \in \Omega$.
\end{itemize} 

If in the definition of \ainc{p} we have $L_p=1$, then we say that $\phi$ satisfies \inc{p}.
Similarly, if $L_q=1$ in the definition of \adec{q}, then we say that $\phi$ satisfies \dec{q}.
By \cite[Section~4.1]{HarH19}, \azero{} can be stated equivalently as the existence of $\beta \in(0,1]$ such that
\[
\phi(x,\beta)
\leq 1
\leq \phi(x,\beta^{-1})
\]
for every $x \in \Omega$.
All the conditions defined above are invariant under the equivalence of $\Phi$-functions.
 
\begin{defn}\label{def:Lphi}
Let $\phi \in \Phiw(\Omega)$ and define the \emph{modular} 
$\rho_\phi$ for $u\in L^0(\Omega)$ by 
\[
\rho_\phi(u)
:= \int_\Omega \phi(x,|u(x)|) \,dx.
\]
The \emph{generalized Orlicz space}, also called Musielak--Orlicz space, is defined as the set 
\[
L^\phi(\Omega)
:= \left\lbrace u \in L^0(\Omega) \colon \lim_{\lambda \to 0^+}\varrho_\phi(\lambda u) = 0 \right\rbrace
\]
equipped with the Luxemburg quasi-norm 
\[
\|u\|_{L^\phi(\Omega)}
:= \inf \left\lbrace \lambda>0 \colon \rho_\phi\left(\frac{u}{\lambda}  \right) \leq 1 \right\rbrace.
\]
If the set is clear from the context we abbreviate $\|u\|_{L^\phi(\Omega)}$ by $\|u\|_{\phi}$.
\end{defn}

A multi-index $\alpha$ is a finite sequence $(a_i)_{i=1}^n$, where every $a_i$ is a non-negative integer.
The length of a multi-index is defined as $|\alpha|:=\sum_{i=1}^n a_i$.
With multi-indices, we can write the partial derivatives of $u$ compactly as
\[
\partial_\alpha u
:= \partial_1^{a_1}\partial_2^{a_2}\cdots\partial_n^{a_n} u,
\]
where $\partial_i := \frac{\partial}{\partial x_i}$ denotes the $i$:th partial derivative.
If $|\alpha| = 0$, then $\partial_\alpha u = u$.

\begin{defn}
A function $u \in L^\phi(\Omega)$ belongs to the
\emph{Orlicz--Sobolev space $W^{k,\phi}(\Omega)$}, if for every multi-index $\alpha$ with $|\alpha| \leq k$ the weak partial derivative $\partial_\alpha u$ exists and belongs to the space $L^{\phi}(\Omega)$.
For $u \in W^{k,\phi}(\Omega)$, we define the quasi-norm
\[
\|u\|_{W^{k,\phi}(\Omega)}
:= \sum_{|\alpha| \leq k} \|\partial_\alpha u\|_{L^\phi(\Omega)}.
\]
If $\Omega$ is clear from the context, we abbreviate $\| u \|_{W^{k,\phi}(\Omega)}$ by $\| u \|_{k,\phi}$.
\end{defn}

%%%%%%%%%%%%%%%%%%%%%%%%%%%%%%%%%%%%%%%%%%%%%%%%%%%%%%%%%%%%%%%%%%%%%%%%%%%%%%%%%%%%%%%%%%%%%%%%%%%%%%%%%%%%%%%%%%%%%%%%%%%%%%

\section{Extension of \texorpdfstring{$\phi$}{phi}}\label{sec:OrliczExtension}

In this section, we study the extension of $\phi \in \Phiw(\Omega)$.
A function $\psi \in \Phiw(\R^n)$ is an extension of $\phi$, if $\psi|_\Omega = \phi$.
We are only interested in extensions, that preserve some of the properties of the original function.
More precisely, we want $\psi$ to satisfy \azero{}, \aone{}, \atwo{} and \adec{q}, if $\phi$ satisfies \azero{}, \aone{}, \atwo{} and \adec{q}.
By \cite[Theorem~3.5]{HarH19}, $\phi$ has an extension satisfying \azero{}, \aone{} and \atwo{}, if and only if $\phi$ satisfies \azero{}, \aomega{} and \atwo{}, where \aomega{} is the following condition:

\begin{itemize}
\item[(A1)$_\Omega$]\label{def:AOmega}
There exists $\beta \in (0,1)$ such that
\[
\beta^{|x-y|t^{1/n}+1}\phi^{-1}(x,t)
\leq \phi^{-1}(y,t)
\]
for every $x,y \in \Omega$ and every $t \geq 1$.
\end{itemize}

By \cite[Lemma~3.3]{HarH19}, \aomega{} is equivalent with \aone{} if $\Omega$ is quasi-convex, and by \cite[Theorem~3.5]{HarH19} the extension satisfies \adec{q} if the original function satisfies it.
Thus, if $\Omega$ is quasi-convex and $\phi$ satisfies \azero{}, \aone{}, \atwo{} and \adec{q} then $\phi$ has an extension satisfying \azero{}, \aone{}, \atwo{} and \adec{q}.
This result however is not sufficient for our purposes, since in our setting $\Omega$ is only an $(\epsilon,\delta)$ domain.

\begin{defn}[from \cite{Jon81}, p.~73]
We say that $\Omega$ is an $(\epsilon,\delta)$-domain, if
for every $x,y \in \Omega$ with $|x-y| < \delta$ there exists a rectifiable curve $\gamma$ such that $\gamma$ lies in $\Omega$ and connects $x$ and $y$ and
\[
\ell(\gamma)
\leq \frac{|x-y|}{\epsilon}
\quad\text{and}\quad
\dist(z,\partial\Omega)
> \frac{\epsilon|x-z||y-z|}{|x-y|}
\]
for every point $z$ lying on $\gamma$.
Let $\Omega_i$ be the connected components of $\Omega$.
We define the radius of $\Omega$ by
\[
\rad(\Omega)
:= \inf_{i}\inf_{x \in \Omega_i}\sup_{y \in \Omega_i}|x-y|.
\]
\end{defn}

To find an extension of $\phi$ in our setting, we use the following approach:
First, we first prove that \aomega{} and \aone{} are equivalent in $(\epsilon,\delta)$-domains.
Then we use \cite[Theorem~3.5]{HarH19} to get an extension of $\phi$ with the desired properties.
We begin by proving that \aomega{} implies \aone{} without any extra assumptions.
The proof is included in the proof of \cite[Lemma~3.3]{HarH19}, but since it is quite simple, we include it here for the sake of completeness.
We use $\omega_n$ to denote the measure of a unit ball.

\begin{lem}\label{lem:AOmegaA1}
Suppose that $\Omega$ is open and $\phi \in \Phiw(\Omega)$ satisfies \aomega{}.
Then $\phi$ satisfies \aone{}.
\end{lem}

\begin{proof}
Suppose first that $\phi$ satisfies \aomega{} with constant $\beta$.
Let $B$ be an arbitrary ball with $|B| \leq 1$.
Let $x,y \in B\cap\Omega$ and $t \in [1,|B|^{-1}]$ be arbitrary.
Let $R$ be the radius of $B$.
Then
\[
|x-y|t^{1/n}
\leq \frac{2R}{|B|^{1/n}}
= \frac{2}{\omega_n^{1/n}}.
\]
Using the above inequality, the fact that $\beta < 1$, and \aomega{}, we find that
\[
\beta^{2\omega_n^{-1/n}+1}\phi^{-1}(x,t)
\leq \beta^{|x-y|t^{1/n}+1}\phi^{-1}(x,t)
\leq \phi^{-1}(y,t).
\]
Thus $\phi$ satisfies \aone{} with constant $\beta^{2\omega_n^{-1/n}+1}$.
\end{proof}

Then we want to prove that \aone{} implies \aomega{} in $(\epsilon,\delta)$-domains, if $\phi$ satisfies \azero{} and \adec{q}.
Instead of $(\epsilon,\delta)$-domains, we can actually formulate the result for somewhat more general domains.

\begin{defn}
Let $\Omega \subset \R^n$, $K \geq 1$ and $\delta > 0$.
We say that $\Omega$ is $(K,\delta)$-quasi-convex, if whenever $x,y \in \Omega$ are such that $|x-y| < \delta$, then there exists a rectifiable curve $\gamma$ connecting $x$ and $y$ such that $\ell(\gamma) \leq K|x-y|$.
\end{defn}

Note that every $(\epsilon,\delta)$-domain is $(\epsilon^{-1},\delta)$-quasi-convex.
We also need the following elementary lemma.

\begin{lem}\label{lem:geqDeltaApu}
Let $\delta > 0$ and $C, q \geq 1$ and $\beta_0\in(0,1)$.
Then there exists $\beta' \in (0,1)$ such that
\[
\beta^{\delta t^{1/n}+1}C^2t
< t^{1/q} 
\]
whenever $t>\beta_0^{-1}$ and $\beta \in (0,\beta')$.
The number $\beta'$ depends on $n$, $\delta$, $C$, $q$ and $\beta_0$.
\end{lem}

\begin{proof}
We will first prove that there exists $\beta' \in (0,1)$ such that
\begin{equation}\label{eq:apulemma}
\frac{\delta t^{1/n}+1}{\ln t}
\geq \frac{1-q}{q\ln\beta}-\frac{2\ln C}{\ln t\ln\beta}
\end{equation}
for $t \in (\beta_0^{-1},\infty)$ and $\beta \in (0,\beta')$.
Since $t \mapsto \frac{\delta t^{1/n}+1}{\ln t}$ is continuous on $(1,\infty)$, and
\[
\lim_{t\to1}\frac{\delta t^{1/n}+1}{\ln t}
=\infty
=\lim_{t\to\infty}\frac{\delta t^{1/n}+1}{\ln t},
\]
there must exist $t_0 \in (1,\infty)$ such that
\[
\frac{\delta t_0^{1/n}+1}{\ln t_0}
= \min_{t \in (1,\infty)} \frac{\delta t^{1/n}+1}{\ln t}
=: m.
\]
Note that $m > 0$, since $t_0 \in (1,\infty)$.
Since 
\[
\lim_{\beta \to 0}\left(\frac{1-q}{q\ln\beta}-\frac{2\ln C}{\ln \beta_0^{-1}\ln\beta}\right)
=0,
\]
there exists $\beta' \in (0,1)$ such that
\[
m
> \frac{1-q}{q\ln\beta}-\frac{2\ln C}{\ln \beta_0^{-1}\ln\beta},
\]
whenever $\beta \in (0,\beta')$.
If $t>\beta_0^{-1}$ and $\beta \in (0,1)$, then
\[
\frac{1-q}{q\ln\beta}-\frac{2\ln C}{\ln \beta_0^{-1}\ln\beta}
> \frac{1-q}{q\ln\beta}-\frac{2\ln C}{\ln t\ln\beta},
\]
since $\beta_0^{-1} > 1$ and $C \geq 1$.
Thus, if $t>\beta_0^{-1}$ and $\beta \in (0,\beta')$, then
\[
\frac{\delta t^{1/n}+1}{\ln t}
\geq m
> \frac{1-q}{q\ln\beta}-\frac{2\ln C}{\ln t\ln\beta}
\]
and \eqref{eq:apulemma} is proved.
Note that $\beta'$ depends only on $m$, $C$, $q$ and $\beta_0$, and $m$ depends only on $n$ and $\delta$.

Multiplying both sides of \eqref{eq:apulemma} by the negative number $\ln t\ln\beta$, and rearranging, we get
\[
(\delta t^{1/n}+1)\ln\beta+2\ln C+\ln t
< \frac{\ln t}{q},
\]
which is equivalent to
\[
\ln (\beta^{\delta t^{1/n}+1}C^2t)
< \ln t^{1/q}. 
\]
Since $\ln$ is an increasing function, it follows that $\beta^{\delta t^{1/n}+1}C^2t < t^{1/q}$.
\end{proof}

\begin{lem}\label{lem:A1AOmega}
Suppose that $\Omega$ is $(K,\delta)$-quasi-convex and that $\phi \in \Phiw(\Omega)$ satisfies \azero{} and \adec{q} with $q \geq 1$.
Then $\phi$ satisfies \aomega{} if it satisfies \aone{}.
\end{lem}

\begin{proof}
Let $\beta_0$ and $\beta_1$ be the constants from \azero{} and \aone{}, respectively, and let $L_q$ be the constant from \adec{q}.
Fix arbitrary points $x,y \in \Omega$, and arbitrary $t \geq 1$.
We must show that the inequality
\begin{equation}\label{eq:AOmega}
\beta^{|x-y|t^{1/n}+1}\phi^{-1}(x,t)
\leq \phi^{-1}(y,t)
\end{equation}
holds for some $\beta \in (0,1)$.
Further $\beta$ must not depend $x$, $y$ and $t$.
We break the proof into several cases.

Suppose first that $x=y$.
Then
\[
\beta^{|x-y|t^{1/n}+1}\phi^{-1}(x,t)
=\beta \phi^{-1}(y,t)
< \phi^{-1}(y,t)
\]
for every $\beta \in (0,1)$.
Thus, in this case, \eqref{eq:AOmega} is true for every $\beta \in (0,1)$.

We then assume that $x \neq y$ and $|x-y| < \delta$.
Since $\Omega$ is a $(K,\delta)$-quasi-convex, there exists a rectifiable curve $\gamma \colon [0,\ell(\gamma)] \to \Omega$, parametrized by arc-length, such that $\gamma$ connects $x$ and $y$ and $\ell(\gamma) \leq K|x-y|$.
Let
\[
M
:= \max \left\lbrace K,\,\omega_n^{-1} \right\rbrace
\]
and let $r := (\omega_n t)^{-1/n}$.
Let $N$ be the unique positive integer such that
\[
\frac{M|x-y|}{r}
< N
\leq \frac{M|x-y|}{r}+1,
\]
and let $x_i := \gamma(i\ell(\gamma)/N)$ for $i = 0,1,\dots,N$.
Since $\gamma$ is parametrized by arc-length, we find that
\[
|x_i-x_{i-1}|
\leq \frac{\ell(\gamma)}{N}
< \frac{r\ell(\gamma)}{M|x-y|}
\leq r\frac{\ell(\gamma)}{K|x-y|}
\leq r.
\]
Thus, for every $i \in \lbrace 1,2,\dots,N \rbrace$ we can find a ball $B_i$ with radius $r$ such that $x_{i-1},x_i \in B_i$.
Since $|B_i|^{-1} = \omega_n^{-1} r^{-n} = t$, \aone{} implies that
\[
\beta_1 \phi^{-1}(x_{i-1},t)
\leq \phi^{-1}(x_i,t),
\]
and we find that
\[
\beta_1^N \phi^{-1}(x,t)
= \beta_1^N \phi^{-1}(x_0,t)
\leq \phi^{-1}(x_N,t)
= \phi^{-1}(y,t).
\]
Using the facts that $\beta_1 < 1\leq M\omega_n^{1/n}$ and $N \leq M|x-y|r^{-1}+1$, we get
\[
(\beta_1^{M\omega_n^{1/n}})^{|x-y|t^{1/n}+1}
\leq \beta_1^{M|x-y|(\omega_n t)^{1/n}+1}
= \beta_1^{M|x-y|r^{-1}+1}
\leq \beta_1^N.
\]
Thus for $\beta\in(0,\beta_1^{2M\omega_n^{1/n}})$ we have
\[
\beta^{|x-y|t^{1/n}+1}\phi^{-1}(x,t)
< \beta_1^N\phi^{-1}(x,t)
\leq \phi^{-1}(y,t).
\]
Hence \eqref{eq:AOmega} is true in the case $|x-y|<\delta$, if $\beta\in(0,\beta_1^{2M\omega_n^{1/n}})$.

Suppose then that $|x-y| \geq \delta$.
Since $\phi$ satisfies \adec{q} by assumption and \ainc{1} by the definition of $\Phiw(\Omega)$, \cite[Proposition~2.3.7]{HarH19b} implies that $\phi^{-1}$ satisfies \adec{1} and \ainc{\frac 1q}.
Thus there exists $L \geq 1$ such that
\begin{equation}\label{eq:aIncADec}
\frac{\phi^{-1}(x,t)}{t}
\leq L\frac{\phi^{-1}(x,1)}{1}
\quad\text{and}\quad
\frac{\phi^{-1}(y,1)}{1}
\leq L\frac{\phi^{-1}(y,t)}{t^{1/q}}.
\end{equation}
Note that $L$ does not depend on $x$, $y$ or $t$.
By \azero{}, we have
\begin{equation}\label{eq:aZeroXY}
\beta_0
\leq \phi^{-1}(x,1)
\leq \frac{1}{\beta_0}
\quad\text{and}{\quad}
\beta_0
\leq \phi^{-1}(y,1)
\leq \frac{1}{\beta_0}.
\end{equation}
We now consider separately the cases $t \in [1,\beta_0^{-1}]$ and $t > \beta_0^{-1}$.

Suppose first that $t \in [1,\beta_0^{-1}]$.
If $\beta \in (0,L^{-1}\beta_0^3)$, then by the fact that $\beta < 1$, \eqref{eq:aIncADec}, \eqref{eq:aZeroXY} and  the fact that $t \mapsto \phi^{-1}(y,t)$ in non-decreasing, we find that
\[
\beta^{|x-y|t^{1/n}+1}\phi^{-1}(x,t)
< \beta \phi^{-1}(x,t)
< \beta Lt\phi^{-1}(x,1)
\leq \beta_0
\leq \phi^{-1}(y,1)
\leq \phi^{-1}(y,t).
\]
Thus \eqref{eq:AOmega} holds in this case for $\beta \in (0,L^{-1}\beta_0^3)$.

Suppose then that $t > \beta_0^{-1}$.
Denote $C:=L\beta_0^{-1}>1$.
Let  $\beta' \in (0,1)$ be the constant given by Lemma \ref{lem:geqDeltaApu}.
We note that $\beta'$ does not depend on $x$, $y$ or $t$.
If $\beta \in (0,\beta')$, then
\[
\beta^{|x-y|t^{1/n}+1}
\leq \beta^{\delta t^{1/n}+1}
\]
since $\delta \leq |x-y|$.
Thus, by \eqref{eq:aIncADec}, \eqref{eq:aZeroXY} and Lemma \ref{lem:geqDeltaApu}, we have
\[
\beta^{|x-y|t^{1/n}+1}\phi^{-1}(x,t)
\leq \beta^{\delta t^{1/n}+1}Lt\phi^{-1}(x,1)
\leq \beta^{\delta t^{1/n}+1}Ct
< C^{-1}t^{1/q}
\leq \phi^{-1}(y,t).
\]
Thus \eqref{eq:AOmega} holds in this case for $\beta \in (0,\beta')$.

Fix some $\beta$ such that
\[
0
< \beta
< \min \left\lbrace \beta_1^{2M\omega_n^{1/n}},\,L^{-1}\beta_0^3,\, \beta' \right\rbrace.
\]
Then \eqref{eq:AOmega} is satisfied.
Since none of the terms $\beta_1^{2M\omega_n^{1/n}}$, $L^{-1}\beta_0^3$ and $\beta'$ depend on $x$, $y$ or $t$, the proof is complete.
\end{proof}

By \cite[Theorem~3.5]{HarH19} and Lemma \ref{lem:A1AOmega}, we get the following corollary:

\begin{cor}\label{cor:phiExtension}
Suppose that $\Omega$ is $(K,\delta)$-quasi-convex, and $\phi \in \Phiw(\Omega)$ satisfies \azero{}, \aone{}, \atwo{} and \adec{q} with $q \geq 1$.
Then $\phi$ has an extension $\psi \in \Phiw(\R^n)$ satisfying \azero{}, \aone{}, \atwo{} and \adec{q}.
\end{cor}

In \cite{HarH19}, the question whether \aone{} and \aomega{} are equivalent in general, was left open.
The next example answers this question by giving an example of a function that satisfies \aone{} but not \aomega{}.
Further, the set in the example is $(\sqrt{2},1)$-quasi-convex, which shows that Lemma \ref{lem:A1AOmega} is not necessarily true without the assumptions \azero{} and \adec{q}.

\begin{eg}
Define the sets $U,V,W \subset \R^2$ by
\[
U
:=\lbrace (x,y) \colon -3<x<-1\text{ and }y>-1\rbrace
,\quad
V
:=\lbrace (x,y) \colon 1<x<3\text{ and }y>-1\rbrace
\]
and
\[
W
:=\lbrace (x,y) \colon -1\leq x\leq 1\text{ and }-1<y<0\rbrace
\]
Note that $\Omega := U\cup V\cup W$ is open and connected.
A simple geometric argument shows that $\Omega$ is $(\sqrt{2},1)$-quasi-convex.

We then define $\phi\in \Phiw(\Omega)$ by
\[
\phi(x,y,t)
:= \begin{cases}
t & \text{if } x \leq 1 \text{ or } y \leq 1, \\
t/y & \text{if } x > 1 \text{ and } y > 1.
\end{cases}
\]
Note that $\phi$ satisfies \dec{1} but not \azero{}.
Also note that $\phi(x,y,t)\neq t$ implies that $(x,y)\in V$.
The inverse is given by
\[
\phi^{-1}(x,y,t)
:= \begin{cases}
t & \text{if } x\leq 1\text{ or }y\leq1, \\
yt & \text{if } x>1\text{ and }y>1.
\end{cases}
\]

Let us show that $\phi$ satisfies \aone{} with constant $\beta = \frac{1}{2}$.
Let $B$ be a ball with $|B| \leq 1$, and let $(x_1,y_1),(x_2,y_2) \in B\cap\Omega$ be arbitrary.
Suppose first that $(x_1,y_1) \notin V$.
Then $x_1 < -1$ or $y_1<0$.
Since the radius of $B$ is less than one, it follows that $x_2 < 0$ or $y_2 < 1$.
Thus
\[
\frac{1}{2}\phi^{-1}(x_1,y_1,t)
=\frac{t}{2}
< t
= \phi^{-1}(x_2,y_2,t).
\]
Similarly we see that $\frac{1}{2}\phi^{-1}(x_1,y_1,t) < \phi^{-1}(x_2,y_2,t)$, if  $(x_2,y_2) \notin V$.

Suppose then that both points belong to $V$.
If $y_1\leq 2$, then
\[
\frac{1}{2}\phi^{-1}(x_1,y_1,t)
\leq t
\leq \phi^{-1}(x_2,y_2,t).
\]
If $y_1 > 2$, then, $2y_2 > 2(y_1-1)>y_1$, and we have
\[
\frac{1}{2}\phi^{-1}(x_1,y_1,t)
=\frac{y_1}{2}t
\leq y_2t
= \phi^{-1}(x_2,y_2).
\]

Let us then show that $\phi$ does not satisfy \aomega{}.
Suppose on the contrary that \aomega{} holds with constant $\beta$.
We have $(-2,y),(2,y)\in\Omega$ for every $y>1$.
Thus \aomega{} implies that
\[
\beta^5 y
=\beta^{|(-2,y)-(2,y)|1^{1/2}+1}\phi^{-1}(2,y,1) 
\leq \phi^{-1}(-2,y,1)
= 1
\]
for every $y>1$.
But this contradicts the fact that $\beta^5y\to\infty$ as $y\to\infty$.
Hence $\phi$ cannot satisfy \aomega.
\end{eg}

%%%%%%%%%%%%%%%%%%%%%%%%%%%%%%%%%%%%%%%%%%%%%%%%%%%%%%%%%%%%%%%%%%%%%%%%%%%%%%%%%%%%%%%%%%%%%%%%%%%%%%%%%%%%%%%%%%%%%%%%%%%%%%

\section{Extension of Sobolev functions}\label{sec:SobolevExtension}

In this section we are going to study the existence of an extension operator $\Lambda: W^{k,\phi}(\Omega) \to W^{k,\psi}(\R^n)$, where $\psi$ is the extension of $\phi$ given by Corollary \ref{cor:phiExtension}.
We will assume that $\Omega$ is $(\epsilon,\delta)$-domain, and that $\phi$ satisfies \azero{}, \aone{}, \atwo{} and \adec{q}.
We will need weighted Sobolev spaces in the process.

\begin{defn}
A weight is a positive locally integrable function $w : \R^n \to (0,\infty]$.
We say that $w$ belongs to Muckenhoupt class $A_1$, if there exists $C \in (0,\infty)$ such that
\[
\frac{1}{|Q|}\int_Q w \,dx
\leq C\essinf_{x \in Q} w(x)
\]
for every cube $Q \subset \R^n$.
A function $u \in L^0(\Omega)$ belongs to the weighted Lebesgue space $L_w^1(\Omega)$, if the norm
\[
\|u\|_{L_w^1(\Omega)}
:= \int_\Omega |u|w\,dx
\]
is finite.
A function $u \in L_w^1(\Omega)$ belongs to the weighted Sobolev space $W_w^{k,1}(\Omega)$, if for every multi-index $\alpha$ with $|\alpha| \leq k$ the weak partial derivative $\partial_\alpha u$ exists and belongs to the space $L_w^{1}(\Omega)$.
For $u \in W_w^{k,1}(\Omega)$, we define the norm
\[
\|u\|_{W_w^{1,k}(\Omega)}
:= \sum_{|\alpha| \leq k} \|\partial_\alpha u\|_{L_w^1(\Omega)}.
\]
\end{defn}

We note that constant weights belong to $A_1$.
The next theorem is the main result of this paper.

\begin{thm}\label{thm:extension}
Let $\Omega$ be an $(\epsilon,\delta)$-domain with $\rad(\Omega)>0$.
Suppose that $\phi \in \Phiw(\Omega)$ satisfies \azero{}, \aone{}, \atwo{} and \adec{q} with $q \geq 1$.
Let $\psi \in \Phiw(\R^n)$ be the extension of $\phi$ given by Corollary \ref{cor:phiExtension}.
Then there exists an operator $\Lambda \colon W^{k,\phi}(\Omega)\to W^{k,\psi}(\R^n)$ and a constant $C$ such that
\[
\|\Lambda u\|_{W^{k,\psi}(\R^n)}
\leq C \|u\|_{W^{k,\phi}(\Omega)},
\]
for every $u \in W^{k,\phi}(\Omega)$.
\end{thm}

\begin{proof}
Let $w$ be an arbitrary weight in the Muckenhoupt class $A_1$.
By \cite[Theorem~1.1]{Chu92}, there exists an extension operator $\Lambda_0: W^{k,1}(\Omega,w) \to W^{k,1}(\R^n,w)$.
Thus
\begin{equation}\label{eq:weightedExtension}
\|\Lambda_0 u\|_{W_w^{k,1}(\R^n)}
\leq C([w]_{A_1})\|u\|_{W_w^{k,1}(\Omega)}
\end{equation}
for all $u \in W_w^{k,1}(\Omega)$.
The constant $C([w]_{A_1})$ depends on $w$, but the dependence is only through $[w]_{A_1}$.
The proof of \cite[Theorem~1.1]{Chu92} shows that $\Lambda_0$ is linear and does not depend on $w$.

Suppose that $u \in W^{k,\phi}(\Omega)$ and that $\supp u$ is bounded.
Then there exists an open ball $B \subset \R^n$ such that $\supp u\subset B$.
Since $|B\cap\Omega|<\infty$, it follows from \cite[Lemma~6.1.6]{HarH19} that $u\in W^{k,1}(B\cap\Omega)$.
Thus $\|u\|_{W^{k,1}(\Omega) }= \|u\|_{W^{k,1}(B\cap\Omega)}<\infty$, and we see that $u\in W^{k,1}(\Omega)$.
It follows that $\Lambda_0 u$ is defined and belongs to $W^{k,1}(\R^n)$.
We define $\Lambda u := \Lambda_0 u$.

In the following $\alpha$ and $\beta$ denote multi-indices.
Define $g \colon \R^n\to[-\infty,\infty]$ by
\[\begin{cases}
g(x)
:= \sum_{|\alpha| \leq k}|\partial_\alpha u|
& \text{if } x \in \Omega, \\
g(x)
:= 0
& \text{if } x \notin \Omega.
\end{cases}\]
If $w \in A_1$, then by \eqref{eq:weightedExtension} we have
\[
\|\partial_\beta(\Lambda u)\|_{L_w^1(\R^n)}
\leq \|\Lambda u\|_{W_w^{k,1}(\R^n)}
\leq C([w]_{A_1})\|u\|_{W_w^{k,1}(\Omega)}.
\]
Since
\[
\|u\|_{W_w^{k,1}(\Omega)}
= \sum_{|\alpha| \leq k} \|\partial_\alpha u\|_{L_w^1(\Omega)}
= \sum_{|\alpha| \leq k} \int_\Omega |\partial_\alpha u|w \,dx
= \int_{\R^n} gw \,dx
= \|g\|_{L_w^1(\R^n)},
\]
we have the inequality
\[
\|\partial_\beta(\Lambda u)\|_{L_w^1(\R^n)}
\leq \|g\|_{L_w^1(\R^n)}.
\]
It then follows from \cite[Corollary~5.3.4]{HarH19b} that there exists a constant $C_0$, such that
\[
\|\partial_\beta(\Lambda u)\|_{L^\psi(\R^n)}
\leq C_0\|g\|_{L^\psi(\R^n)}.
\]
Since $\psi|_\Omega = \phi$ and $g=0$ outside $\Omega$, $C_0\|g\|_{L^\psi(\R^n)} = C_0\|g\|_{L^\phi(\Omega)}$.
Denote	by $c_{n,k}$ is the number of multi-indices with $|\beta| \leq k$.
Using the above estimates and \cite[Corollary~3.2.5]{HarH19b} we find that
\[
\sum_{|\beta| \leq k}\|\partial_\beta(\Lambda u)\|_{L^\psi(\R^n)}
\leq C_0 \sum_{|\beta| \leq k}\|g\|_{L^\phi(\Omega)}
= c_{n,k}C_0\|g\|_{\phi}
\leq c_{n,k}C_0L\sum_{|\alpha|\leq k}\|\partial_\alpha u\|_{\phi},
\]
where $L \geq 1$ is the constant from \cite[Corollary~3.2.5]{HarH19b}.
Hence
\begin{equation}\label{eq:boundedSuppExtension}
\|\Lambda u\|_{W^{k,\psi}(\R^n)}
\leq C\|u \|_{W^{k,\phi}(\Omega)}
\end{equation} with $C :=  c_{n,k}C_0L$.

Let then $u \in W^{k,\phi}(\Omega)$ be such that $\supp u$ is unbounded.
By \cite[Lemma~6.4.1]{HarH19b}, the set of functions in $W^{1,\phi}(\Omega)$ with bounded support is dense in $W^{1,\phi}(\Omega)$.
A similar proof works also for $k > 1$ (see also \cite[Theorem~9.1.2]{DieHHR11}).
Thus there exist a sequence of functions $u_i \in W^{k,\phi}(\Omega)$, such that $\supp u_i$ is bounded and $u_i\to u$ in $W^{k,\phi}(\Omega)$.
Further, we see from the proof of \cite[Lemma~6.4.1]{HarH19b} that we can choose $u_i$ to be such that $u_i|_{B(0,i)\cap \Omega} = u|_{B(0,i) \cap \Omega}$ for every $i$.

By the linearity of $\Lambda_0$ and \eqref{eq:boundedSuppExtension}, we have
\[
\|\Lambda u_i-\Lambda u_j\|_{W^{k,\psi}(\R^n)}
= \|\Lambda(u_i-u_j)\|_{W^{k,\psi}(\R^n)}
\leq C\|u_i-u_j \|_{W^{k,\phi}(\Omega)}.
\]
Since $(u_i)$ is a Cauchy sequence in $W^{k,\phi}(\Omega)$, it follows that $\Lambda u_i$ is a Cauchy sequence in $W^{k,\psi}(\R^n)$.
Since $W^{k,\psi}(\R^n)$ is a Banach space, there exists $v\in W^{k,\psi}(\R^n)$ such that $\Lambda u_i\to v$.
Note that $v|_\Omega = u$, since
\[
\Lambda u_i|_{B(0,i)\cap \Omega} = u_i|_{B(0,i)\cap \Omega}
= u|_{B(0,i)\cap \Omega}
\]
for every $i$.
Thus, we can extend $u$ by defining $\Lambda u:= v$.
Then $\Lambda u_i\to\Lambda u$ in  $W^{k,\psi}(\R^n)$.
Let $\eta > 1$ be arbitrary.
Since $u_i\to u$ and $\Lambda u_i\to\Lambda u$, for large enough $i$ we have
\[
\|\Lambda u\|_{W^{k,\psi}(\R^n)}
\leq \eta\|\Lambda u_i\|_{W^{k,\psi}(\R^n)}
\leq \eta C\|u_i\|_{W^{k,\phi}(\Omega)}
\leq \eta^2 C\|u\|_{W^{k,\phi}(\Omega)}.
\]
Letting $\eta \to 1$ completes the proof.
\end{proof}

\begin{rem}
The function $v$ in the previous proof is independent of the sequence $(u_i)$.
Let $(w_j)$ be another sequence of functions with compact support approaching $u$ in $W^{k,\phi}(\Omega)$.
Let $w \in W^{k,\phi}(\R^n)$ be the function that $\Lambda w_i$ approaches.
Let $\epsilon > 0$.
By the the quasi-triangle inequality,
\[
\|v - w\|_{W^{k,\psi}(\R^n)}
\lesssim \|v - \Lambda u_i\|_{W^{k,\psi}(\R^n)} + \|\Lambda u_i - \Lambda w_j\|_{W^{k,\psi}(\R^n)} + \|\Lambda w_j - w \|_{W^{k,\psi}(\R^n)}.
\]
If $i$ and $j$ are large enough, then
\[
\|v - w\|_{W^{k,\psi}(\R^n)}
\lesssim \|\Lambda u_i - \Lambda w_j\|_{W^{k,\psi}(\R^n)} + \epsilon.
\]
By the linearity of $\Lambda_0$, \eqref{eq:boundedSuppExtension} and the quasi-triangle inequality,
\[
\begin{aligned}
& \|\Lambda u_i - \Lambda w_j\|_{W^{k,\psi}(\R^n)}
= \|\Lambda(u_i-w_j)\|_{W^{k,\psi}(\R^n)}
\leq C\|u_i - w_j \|_{W^{k,\phi}(\Omega)} \\
& \lesssim C\|u_i - u \|_{W^{k,\phi}(\Omega)} +  C\|u - w_j \|_{W^{k,\phi}(\Omega)}
< \epsilon
\end{aligned}
\]
for sufficiently large $i$ an $j$.
Thus $\|v - w\|_{W^{k,\psi}(\R^n)} \lesssim \epsilon$.
Since $\epsilon > 0$ was arbitrary, $v = w$.
\end{rem}

If $\psi \in \Phiw(\R^n)$ satisfies \azero, \aone{}, \atwo{} and \adec{q}, then  $C_0^\infty(\R^n)$ is dense in $W^{1,\psi}(\R^n)$ by \cite[Theorem~6.4.4]{HarH19b}.
Combining this with Theorem \ref{thm:extension}, we get the following corollary:

\begin{cor}
Let $\Omega$ be an $(\epsilon,\delta)$-domain and suppose that $\phi \in \Phiw(\Omega)$ satisfies \azero, \aone{}, \atwo{} and \adec{q} with $q \geq 1$.
Then $C_0^\infty(\overline{\Omega})$ is dense in $W^{1,\phi}(\Omega)$.
\end{cor}

%%%%%%%%%%%%%%%%%%%%%%%%%%%%%%%%%%%%%%%%%%%%%%%%%%%%%%%%%%%%%%%%%%

\end{document}